\documentclass[12pt]{article}
\usepackage{hyperref}

\usepackage[T1]{fontenc}
\usepackage[utf8]{inputenc}
\usepackage{authblk}

\usepackage{dsfont}
\usepackage{amssymb,amsthm,amsmath}

\newtheorem{thm}{Theorem}[section]

 \makeatletter

\@addtoreset{lemma}{thm}

\newtheorem{corollary}{Corollary}[section]
 \makeatletter

\@addtoreset{corollary}{thm}

\newtheorem{remark}{Remark}[section]

\begin{document}

\date{}
\title{Asymptotics of Randomly Weighted $u$- and $v$-statistics: Application to Bootstrap  }

\author{Mikl\'{o}s Cs\"{o}rg\H{o}\thanks{mcsorgo@math.carleton.ca} and Masoud M. Nasari\thanks{mmnasari@math.carleton.ca}    \\
\small{School of Mathematics and Statistics of Carleton University}\\\small{ Ottawa, ON, Canada} }

\maketitle
\begin{abstract}
\noindent
This paper is mainly concerned with asymptotic studies  of weighted bootstrap for $u-$ and $v-$statistics. We derive the consistency of the weighted bootstrap $u-$ and $v-$statistics, based on i.i.d. and non i.i.d. observations,   from some more general results which we first establish for sums of randomly weighted arrays of random variables. Some of the results in this paper significantly extend some  well-known results on consistency of $u$-statistics and also consistency of  sums of arrays of random variables. We also employ a new approach to conditioning  to  derive a conditional CLT for weighted bootstrap $u-$ and $v-$statistics,  assuming the same conditions as the classical central limit theorems for regular $u-$ and $v-$statistics.
\end{abstract}



\section{Introduction}
\label{section 1}
The main purpose of this study is to investigate the validity of bootstrap $u$- and $v$-statistics resulting for the so-called $m$-out of-$n$   scheme  of bootstrap. The use of this scheme results in having  a randomly weighted  form of the original $u$, $v$-statistics. This, in turn,  exhibits  an explicit connection between the original $u$, $v$-statistics and their  bootstrap version. Moreover, it clearly identifies the two types of variations involved in a bootstrap problem for statistics of the form of partial sums in general and for $u$, $v$-statistics in particular.  In this exposition, the investigation of bootstrap $u$, $v$-statistics  includes the consistency and convergence in distribution of  these bootstrap statistics.  While the results for investigating the convergence in distribution  for bootstrap is done only for $u$, $v$-statistics, based on i.i.d. observations, when $m$-out of-$n$ scheme is used, the results in this paper on the consistency are provided in a more general way. In fact, the results on the consistency is restricted to neither $u$, $v$-statistics nor to i.i.d. observations or bootstrap. The results  on consistency are dealt with for randomly weighted  arrays of random variables. These results then are used to derive consistency for randomly weighted $u$-statistics as well as for bootstrap $u$, $v$-statistics based on i.i.d. observations and observations with absolute regularity property (cf. Theorems \ref{Corollary u-stat iid} and \ref{Corollary u-stat beta-mixing}, respectively). Some of the results in this work extend  well-known results on the strong law of large numbers for arrays of random variables and also for $u$-statistics (cf. Theorem \ref{revers martingale} and Remark \ref{Berk's}),  and some of them shed  light on the consistency of bootstrap when the original  sample can stay  finite  (cf. Theorems \ref{SLLN for arrays fixed n} and \ref{slln for u-stat fix n}).
\par
The  conditions assumed for the results in this paper are the same as the ones required  for the classical central limit  theorems (CLTs)  and classical strong and  weak laws of large numbers in the non-weighted case. In other words, there are  no further restrictions imposed  on the observations.
\par
The material in this paper is organized as follows. In Section \ref{section 2} laws of large numbers are provided  for randomly weighted arrays of random variables  (observations).   The presentation of these results is  so that they can be used in establishing the consistency of bootstrap $u$, $v$-statistics which is provided in Section \ref{section 3}. Section \ref{section 4} is devoted to establishing conditional (given the weights) central limit theorems for $u$, $v$-statistics. In Section \ref{section general case}, remarks are made  on the validity of the results for higher dimensional arrays of random variables and $u$, $v$-statistics of order greater than $2$. The proofs are given in Section \ref{section 5}.

\section{Laws of large numbers for randomly weighted arrays of random variables}
\label{section 2}
Consistency of  bootstrap mean as  randomly weighted sums was pioneered by Athreya \cite{Athreya} and followed by S. Cs\"{o}rg\H{o} \cite{S. Csorgo} and Arenal-Guti\'{e}rrez \emph{et al}.   \cite{Arenal et al}. Since then the problem  has received a great deal of  attention from researchers. Of the contributions to the field, we specifically mention the two results, respectively  due Arenal-Guti\'{e}rrez \emph{et al}. \cite{Arenal et al},  and  Rosalsky and Sreehari \cite{rosalsky}, as the former influenced Theorem \ref{slln for u-stat fix n}  and  the latter  motivated    Theorem \ref{SLLN for arrays} of this exposition.
\par
In this paper  we consider  two sequences of possibly double  triangular (with respect to $m$ and $n$) arrays of random variables, $\{X^{(n,m)}_{ij}; \ 1\leq i,j\leq n \}$  and $\{\varepsilon^{(n,m)}_{ij}; \ 1\leq i,j\leq n\}$, $n,m\geq 1$,   which are  defined  on  the same  probability space
$(\Omega_{X,\varepsilon}, \mathfrak{F}_{X,\varepsilon}, P_{X,\varepsilon})$.
 Also, by $(\Omega_{X},\mathfrak{F}_{X},P_{X})$ and
$(\Omega_{\varepsilon},\mathfrak{F}_{\varepsilon},P_{\epsilon})$  we denote the  marginal probability spaces of the  $X^{(n,m)}$s and   $\varepsilon^{(n,m)}$s, respectively. We shall, often,  refer to the $X^{(n,m)}$s as the observations (data) and $\varepsilon^{(n,m)}$s are referred to  as the weights. We shall investigate the large sample behavior of the   randomly weighted sums $\sum_{1\leq i,j\leq n} \varepsilon^{(n,m)}_{ij} X^{(n,m)}_{ij}$, $n,m\geq 1$.   The observations and the weights are so that they can be employed in studying the $m$-out of-$n$ scheme of bootstrap for $u$-statistics which is to be discussed   in Sections \ref{section 3} and  \ref{section 4}.
We note in passing that when the observations and the weights are independent, then   their  joint probability space can of course  be defined as the direct product probability space    $(\Omega_{X}\times \Omega_{\varepsilon},  \mathfrak{F}_{X} \otimes \mathfrak{F}_{\varepsilon}, P_{X}\times P_{\varepsilon})$ of their marginals.
\par
In this section we present some strong and weak consistency   results  for  sums of randomly weighted arrays of random variables.
\par
Except for Theorem \ref{SLLN for arrays} below and  its application to bootstrap $u$-statistics in  Theorem \ref{Marcinkiewicz} of the next section,   the method of conditioning plays an important role in the  establishment of the results in this exposition. More precisely, employing  hierarchical   arguments,  we derive our results   via conditioning  on the weights $\varepsilon^{(n,m)}$ in some stochastic way with respect to $P_{\varepsilon}$. The latter results, in turn, can be  extended to unconditional ones in terms of  the joint probability measure $P_{X,\varepsilon}$. Hence, we let $P_{.|\varepsilon}(.)$ and $E_{.|\varepsilon}(.)$, respectively stand for the  conditional probability and conditional expected value given the weights $\varepsilon^{(n,m)}$.
\par
The following Theorem \ref{SLLN for arrays} is a strong law of large numbers for sums of randomly weighted arrays of random variables.
\begin{thm}\label{SLLN for arrays}
Consider the two possibly (double) triangular arrays $\{X^{(n,m)}_{ij};$ $1\leq i,j \leq n\}$  and $\{\varepsilon^{(n,m)}_{ij};\ 1\leq i,j \leq n \}$, $n,m\geq 1$, of random variables which are defined on the same probability space. Let the sequence of positive integers $m$ be such that  $m=m(n) \rightarrow +\infty$, as $n\rightarrow +\infty$. Also let  $\{c^{(n,m)}_{ij}; \ 1\leq i,j\leq n  \}$ be a possibly (double) triangular array of real numbers and $\{a^{(n,m)}_{ij}; \ 1\leq i,j\leq n  \}$ be a possibly (double) triangular array of positive real numbers such that, for $\delta>0$,
\begin{eqnarray}
&&\hspace{-4.5 cm} (a)\ P_{\varepsilon}( \cup_{1\leq i,j\leq n} | \varepsilon^{(n,m)}_{ij} - c^{(n,m)}_{ij} | > \delta \ a^{(n,m)}_{ij}, \ i.o. (n)\ )=0, \label{eq 1}  \\
&&\hspace{-4.5 cm}  and, as\ m(n),n \rightarrow  +\infty,\nonumber
\\
&&\hspace{-4.5 cm}  (b)\ \sum_{1\leq i,j\leq n} a^{(n,m)}_{ij} |X^{(n,m)}_{ij}|<+\infty \ a.s.-P_X, \nonumber \\
&& \hspace{-4.5 cm} (c)\ \sum_{1\leq i,j\leq n} c^{(n,m)}_{ij} X^{(n,m)}_{ij}<+\infty \ a.s.-P_X. \nonumber
\end{eqnarray}
Then,  (a), (b) and (c), as $n\rightarrow \infty$, imply that
\\ \\
$\sum_{1\leq i,j\leq n} \varepsilon^{(n,m)}_{ij} X^{(n,m)}_{ij}$ converges a.s.-$P_{X,\varepsilon}$ to the same a.s.-$P_X$ limit as that of  $\sum_{1\leq i,j\leq n} c^{(n,m)}_{ij} X^{(n,m)}_{ij}$.
\end{thm}
\noindent
\begin{remark}
In many cases it is  natural  for the numerical  sequence $c_{i,j}^{(n,m)}$ to be taken to be the mean of the  weights, when it exists and is finite. In other words, for each,  $1\leq i,j\leq n $, $c_{i,j}^{(n,m)}=E_{\varepsilon} \varepsilon_{ij}^{(n,m)}$, $n,m \geq 1$.  In this case, in view of the first  Borel-Cantelli Lemma  and  Chebyshev inequality, perhaps the first and  most natural way to investigate (\ref{eq 1}), when the weights have finite second moments,  is  to check if the following  holds true:
$$ \sum_{n=1}^{\infty}\sum_{1\leq i,j\leq n} \frac{Var(\varepsilon_{ij}^{(n,m)})}{(a_{ij}^{(n,m)})^2}<+\infty.$$
This would amount to a generalization of   Remark 1 of  \cite{rosalsky} to the case of arrays of random variables. However, in this exposition, when studying bootstrap $u$- and $v$-statistics in Section \ref{section 3}, we shall use Bernstien's inequality to investigate (\ref{eq 1}), as it gives sharper bounds with less computation when the weights are products of  multinomially distributed random variables.
\end{remark}

Theorem \ref{SLLN for arrays} above was motivated by, and it generalizes,  Theorem 1 of Rosalsky and Sreehari \cite{rosalsky} to arrays of random  variables so that it can be used in studying the  validity of the so called  $m$-out of-$n$ method of bootstrap $u$-statistics in Section \ref{section 3}. It is obvious that  when dealing with regular triangular arrays of random variables then, on taking  $\varepsilon^{(n,m)}=\varepsilon^{(n)}$ and $X^{(n,m)}=X^{(n)}$, Theorem \ref{SLLN for arrays} is true. It also continues to hold true for non-triangular observations and arrays.

The following Theorem \ref{SLLN for arrays fixed n} assumes that the observations have  finite means  and it concerns the randomly weighted sums of the form
$$\sum_{1\leq i,j \leq n}  \varepsilon^{(n,m)}_{ij} X^{(n,m)}_{ij}$$
of possibly (double) triangular arrays of observations
$\{X^{(n,m)}_{ij}; 1\leq i,j\leq n \}$
and possibly (double) triangular arrays of random weights $\{\varepsilon^{(n,m)}_{ij}; 1\leq i,j\leq n \}$, $n,m\geq 1$. The key difference between Theorem \ref{SLLN for arrays} above and the next result is that it does not necessarily  require that both $n,m\rightarrow +\infty$. In fact it can be true when both, or either one of $m$ and $n$, approach   $+\infty$. As a consequence of this, it  leads to an interesting application to  bootstrap $u$, $v-$statistics, via the $m$-out of-$n$ scheme,   when the number of observations $n$ is  fixed and only the bootstrap sample size $m$ approaches $+\infty$. The latter result is presented in  Theorem \ref{slln for u-stat fix n} of the next section.

\begin{thm}\label{SLLN for arrays fixed n}
With  the positive integers $n,m\geq 1$, let $\{X^{(n,m)}_{ij}; \ 1\leq i,j\leq n \}$ and $\{\varepsilon^{(n,m)}_{ij}; \ 1\leq i,j\leq n \}$, which are independent from each other, be possibly (double) triangular arrays of random variables such that, for  $1\leq i,j\leq n$,  $E_{X}|X^{(n,m)}_{ij}|\leq c^{(n,m)}_{ij}$, where $c^{(n,m)}_{ij}>0$. Consider the following four statements:
\\
(i)  $n$ is fixed and $m\rightarrow+\infty$,
\\
(ii) $m$ is fixed and $n\rightarrow+\infty$,
\\
(iii) $n,m\rightarrow+\infty$ such that $m=m(n)$ is an increasing function of $n$,
\begin{equation*}\label{change to probability 1}
(iv)\ \sum_{1\leq i,j\leq n} \big|\varepsilon^{(n,m)}_{ij}\big| c^{(n,m)}_{ij}=o(1)\ a.s.-P_{\varepsilon}. \qquad \qquad \qquad \qquad \qquad \qquad \qquad \qquad
\end{equation*}
Then,   either (i) and (iv), or (ii) and (iv), or  (iii) and (iv) suffice to have
\begin{equation}\label{change to probability 2}
P_{X|\varepsilon} \big( \big|\sum_{1\leq i,j \leq n}  \varepsilon^{(n,m)}_{ij} X^{(n,m)}_{ij}    \big|>\delta  \big)\rightarrow 0 \ a.s.-P_{\varepsilon},
\end{equation}
for any $\delta>0$.
\\
Also, when  (iv) holds true in probability-$P_{\varepsilon}$ to begin with,  then so is the conclusion (\ref{change to probability 2}).
\end{thm}

\begin{remark}\label{new 1}
Observe that when the observations are non-triangular, then naturally  for the bounds we have  $c_{ij}^{(n,m)}=c_{ij}$. In this case if $ \sup_{1\leq i,j\leq +\infty} c_{ij}<+\infty$, then (iv) of Theorem \ref{SLLN for arrays fixed n} can be replaced by
\begin{equation}\nonumber
\sum_{1\leq i,j\leq n} \big|\varepsilon^{(n,m)}_{ij}\big| =o(1)\ a.s.-P_{\varepsilon}.
\end{equation}
\end{remark}

\begin{remark}\label{used for the corollary}
Since $ P_{X|\varepsilon} \big( \big|\sum_{1\leq i,j \leq n}  \varepsilon^{(n,m)}_{ij} X^{(n,m)}_{ij}    \big|>\delta  \big)  \leq 1$, by virtue of Lemma 1.2 in  S. Cs\"{o}rg\H{o} and Rosalsky \cite{Csorgo Rosalsky},
conclusion (\ref{change to probability 2}) of Theorem \ref{SLLN for arrays fixed n},  implies that
\begin{equation}\nonumber
\sum_{1\leq i,j \leq n}  \varepsilon^{(n,m)}_{ij} X^{(n,m)}_{ij}\rightarrow 0, \ in \ probability-P_{X,\varepsilon}.
\end{equation}

\end{remark}

The following result establishes a strong law of large numbers for sums of  randomly weighted arrays of random variables when neither the observations nor the random random weights are triangular while both are assumed to be symmetric. This result is a generalization of the strong law of large numbers for $u$-statistics due to Berk \cite{Berk} which was proven in view of martingale property of $u$-statistics. Some ideas were borrowed  from the latter paper and adapted accordingly.

\begin{thm}\label{revers martingale}
Let $X_{ij}$ and $\epsilon_{ij}$, $1\leq i\neq j\leq n$, $n\geq 1$,  be \emph{symmetric} arrays of random variables which are defined on the same probability space. If,  for all $i,j$,   $\epsilon_{ij} X_{ij}$ are identically distributed and $E_{X,\epsilon}|\epsilon_{12}X_{12}|<+\infty$, then, as $n\rightarrow \infty$,
\begin{equation*}
S_{n}:=\frac{1}{n(n-1)} \sum_{1\leq i\neq j\leq n} \epsilon_{ij}X_{ij}\rightarrow E_{X,\epsilon}(\epsilon_{12}X_{12})\ a.s.-P_{X,\varepsilon}.
\end{equation*}
Trivially, when  $\epsilon_{ij}$ and $X_{ij}$ are uncorrelated with $E_{\epsilon}| \epsilon_{ij}|<+\infty$,  and $E_{X}| X_{ij}|<+\infty$,    then the  limit above becomes $E_{\epsilon}(\epsilon_{12}) E_{X}(X_{12})$.
\end{thm}
\begin{remark}\label{Berk's}
Theorem \ref{revers martingale}  generalizes Berk's \cite{Berk} strong law of large numbers for regular $u$-statistics based on i.i.d. observations to randomly weighted arrays of random variables that may be taken to be the summands  of a (randomly) weighted  $u$-statistic based on identically distributed but not necessarily independent observations.
\end{remark}
Theorem \ref{revers martingale} also  generalizes Theorem 2 of Etemadi \cite{Etemadi}, on  taking $m=n$ in the latter. The latter result of  Etemadi requires that $\epsilon_{ij} X_{ij}$s  be pairwise independent and  identically
distributed random variables, with respect to $P_{X,\epsilon}$, and that $E_{X,\epsilon}\{|\epsilon_{ij} X_{ij}|\log|\epsilon_{ij} X_{ij}| \}< +\infty$. As it can be readily seen, Theorem \ref{revers martingale}, drops the independence as well as it requires only $E_{X,\epsilon} |\epsilon_{ij} X_{ij}| < +\infty$.

\par
We note in passing that  Theorems \ref{SLLN for arrays}, \ref{SLLN for arrays fixed n} and \ref{revers martingale} assume almost no  conditions  (such as requiring certain kind of dependencies, for example) on the relation between  the weights and the observations.

\section{Consistency of bootstrap $u$- \& $v$-statistics}\label{section 3}
Studying the problem of consistency of the bootstrap usually consists of showing that the deviation of the statistic in hand and its bootstrap version, perhaps multiplied by some normalizing sequence, vanishes  a.s. or in probability with resect to the joint probability of the observations and the bootstrap experiment. In large sample applications of the  bootstrap, one use of such  an  approximation  is  establishing  limiting distribution results for a bootstrap statistic in hand. As an example in our context,  we refer to the consistency   in relation (\ref{eq proof 7}) below, noting that $G_{n,m}$ is a deviation between a  $u$-statistic with the kernel $h(X_1,X_2)-\tilde{h}(X_1)-\tilde{h}(X_2)$ and its weighted bootstrap version.
The need for such consistencies has first become apparent  in  A(6) of Burke and Gombay \cite{Burk and Gombay}. Motivated by this, and in view of our results in the previous section,  in this section we give some results on  the problem of bootstrapped  $u-$ and $v-$statistics when the original  sample  is large.
\par
In small sample theory, the problem of consistency of bootstrap, as explained right above Theorem \ref{SLLN for arrays fixed n},  leads to interesting results, if one can show that by re-sampling repeatedly from a finite sample,  the thus obtained sequence of bootstrap versions of a statistic in hand  converges to the original statistic  a.s. or in probability, with respect to the  joint probability of the observations and the bootstrap experiment.  The results which we developed in Section \ref{section 2}  enable  us to make such conclusions for bootstrapped  $u$- and $v-$statistics when the original sample is finite (cf. Theorems \ref{SLLN for arrays fixed n} and \ref{slln for u-stat fix n}).
\par
In this section, as well as in the next one, we shall study  the so-called  $m$-out of-$n$ scheme of bootstrap which we are now to detail as follows.
\section*{$m$-out of-$n$ scheme of bootstrap}
Draw a sample of size $m\geq 1$  with replacement from the set $\{1,\ldots,n \}$  of the indices of the original observations $X_1,\ldots,X_n$. In this scheme a bootstrap  sub-sample of size $m$, denoted by $X^{*}_{1},\ldots,X^{*}_{m}$, is drawn from the original sample $X_1,\ldots,X_n$, $n\geq 1$,  independently  of the original observations.
\par
 In this section  we  study the consistency of bootstrap $u$- and $v-$statistics via the scheme of $m$-out of-$n$. The obtained results  will be seen to be applicable to i.i.d. and some stationary  observations.
\par
A $u$-statistic of order 2 is defined by
\begin{equation}\label{definition of u-satat}
U_n= \frac{1}{n(n-1)} \sum_{1\leq i\neq j \leq n}  h(X_i,X_j),
\end{equation}
where, $h:\mathds{R}^2 \rightarrow \mathds{R}$ is a measurable function which is symmetric in its arguments,  and is called the kernel.  Many of the well known statistics, such as sample variance, deleted jackknife variance estimator, Fisher's \textit{k}-statistic for estimation of cumulants, Kendal's $\tau$, Gini's mean difference,  are  examples of $u$-statistics.
\par
When investigating  bootstrapped $u$-statistics, it is quite  natural to consider $v$-statistics as well. A $v$-statistic
can be viewed as an extension of a $u$-statistic of the same order,  as it is defined as follows.
\begin{eqnarray}
V_n&=& \frac{1}{n^2} \sum_{1\leq i,j \leq n} h(X_i,X_j).\label{definition of v-satat}\\
&=& \frac{n-1}{n} \ U_n + \frac{1}{n^2}  \sum_{i=1}^{n} h(X_i,X_i), \label{relation u- and v}
\end{eqnarray}
where $U_n$ is the $u$-statistic as in (\ref{definition of u-satat}). The major  difference between the $u$-statistic $U_n$ and the $v$-statistic $V_n$ is that the latter one includes the  diagonal terms $h(X_i,X_i)$, $1\leq i \leq n$, while the former one doesn't.  When the kernel $h$ is so that $h(x,y)=\ell(x-y)$, where $\ell$ is a function, then the corresponding  $u$- and $v$-statistics coincide up to the constant $(n-1)/n$. Even when this is not the case,  in view of (\ref{relation u- and v}), under some regularity  conditions, $U_n$ and $V_n$ will asymptotically  coincide.
\par
The reason that one should consider $v$-statistics when  bootstrapping $u$-statistics using the scheme of $m$-out of-$n$  is that, as a result of re-sampling with replacement, when computing the bootstrap $u$-statistic, for $1\leq s\neq t\leq m$, we may for example have $X^{*}_{s}=X_{1}$ and $X^{*}_{t}=X_1$. This event has  probability $1/n^2$ that vanishes as $n$ increases  to $+\infty$.
\par
A bootstrap $u$-statistic $U^{*}_{n,m}$ and bootstrap $v$-statistic $V^{*}_{n,m}$     based on a bootstrap sub-sample $X^{*}_1,\ldots,X^{*}_{m}$ of $X_1,\ldots,X_n$, resulting from the method of $m$-out of-$n$ of bootstrap,  are respectively  defined  as follows.
\begin{eqnarray*}
U^{*}_{n,m}&:=& \frac{\sum_{{\substack{1\leq s\neq t \leq m \\ X^{*}_{t} and X^{*}_{s}\ are \ distinct  } }} h(X^{*}_s,X^{*}_t)} {m(m-1)}\\
&=&\frac{\sum_{1\leq i\neq j \leq n} w^{(n)}_{i} w^{(n)}_{j} h(X_{i},X_{j})   }{m(m-1)}\nonumber
\end{eqnarray*}
and
\begin{eqnarray}
V^{*}_{n,m}&:=& \frac{\sum_{1\leq s\neq t \leq m} h(X^{*}_s,X^{*}_t)} {m(m-1)}\label{Def of boostrap u-stat}\\
&=& \frac{\sum_{1\leq i, j \leq n} w^{(n)}_{i} w^{(n)}_{j} h(X_{i},X_{j}) }{m(m-1)},\nonumber
\end{eqnarray}
where $w^{(n)}_{i}$ is the number of times the index $i$, $1\leq i \leq n$, is chosen in the scheme of $m$-out of-$n$. In view of our earlier definition of this scheme of bootstrap it is  obvious that  $w_{i}^{(n)}$s are independent from the original observations $X_{1},\ldots,X_n$. Also, it is easy to see that  $w^{(n)}_{i}$, $1\leq i \leq n$, to which we shall refer to  as bootstrap weights or simply  weights, are so that $\sum_{i=1}^{n} w^{(n)}_{i}=m$ and  $E_{w}(w^{(n)}_{i}/m)=1/n$, $1\leq i \leq n$. That is, for each $n\geq 1$, the bootstrap weights have  multinomial distribution with size $m$. In other words,
\begin{equation}\nonumber
(w^{(n)}_{1},\ldots,w^{(n)}_{n})\  \substack{d\\=}\ \ multinomial(m;\frac{1}{n},\ldots,\frac{1}{n}).
\end{equation}
\begin{remark}\label{relation between bootstrap u,v}
To state our results for $U^{*}_{n,m}$ and $V^{*}_{n,m}$, it is important to address the relation between the two bootstrap statistics. Simple calculations show that, rhyming with (\ref{relation u- and v}), we have
\begin{equation}\label{eq relation bootstrap u,v}
V^{*}_{n,m}= U^{*}_{n,m}+ \sum_{i=1}^n \frac{(w_{i}^{(n)})^2}{m(m-1)} h(X_i,X_i).
\end{equation}
\end{remark}
\par
In this section, for a variety of  observations, independent or dependent, we shall show that, as $n,m\rightarrow +\infty$,  and at some point  when only $m\rightarrow +\infty$, the deviation between the bootstrap $u$-statistic $U^{*}_{n,m}$ and the associated original $u-$ and  $v-$statistics $U_n$ and $V_n$ goes to zero a.s.,  or in probability,  with respect to the joint distribution of $P_{X,w}$.
\par
For the following result we consider i.i.d. observations and, as both $m,n\rightarrow +\infty$, we present a Marcinkiewicz type law of large numbers for $U^{*}_{n,m}$ and $V^{*}_{n,m}$  when the kernel $h$ has less than one moment. In fact Theorem \ref{Marcinkiewicz}, below,  can be viewed as a bootstrap version of Theorem 1 of Gin\'{e} and Zinn \cite{Gine Zinn}. This result is a consequence of Theorem \ref{SLLN for arrays}  and  reads as follows.
\begin{thm}\label{Marcinkiewicz}
Let $X_1,\ldots,X_{n}$ be the first $n\geq 1$ terms of an infinite sequence of  i.i.d. random variables. \\
(a) Assume that $E_{X}|h(X_1,X_2)|^{2/d}<+\infty$, where $d>2$. Then, as $m,n\rightarrow +\infty$ such that $m=O(n^{d/(d-2) }\log^{\frac{d+2}{d-2}} n)$,
\begin{equation}\nonumber
m^{-d+2}\ U^{*}_{n,m}\rightarrow 0 \ a.s.-P_{X,w}.
\end{equation}
(b) Assume that $E_{X}|h(X_1,X_2)|^{2/d}<+\infty$ and $E_{X}|h(X_1,X_1)|^{2/d}<+\infty$, where $d>2$. Then, as $m,n\rightarrow +\infty$ such that $m=O(n^{d/(d-2) }\log^{\frac{d+2}{d-2}} n)$,
\begin{equation}\nonumber
m^{-d+2}\ V^{*}_{n,m}\rightarrow 0 \ a.s.-P_{X,w}.
\end{equation}
\end{thm}
\par
When the first moments of   the kernel  $h$ of $U_n$ and  $V_n$, which does not depend on the sample size $n$,   are uniformly bounded in $1\leq i,j <+\infty$,  then the following Theorem \ref{slln for u-stat fix n} is a  direct consequence to Theorem \ref{SLLN for arrays fixed n} and Remark \ref{used for the corollary}.

\begin{thm}\label{slln for u-stat fix n}
Let $X_1, \ldots,X_{n}$, $n\geq 1$, which are the first $n\geq 1$ terms of  an infinite  sequence of random variables, and the kernel $h$ be such that
\begin{equation}\nonumber
\sup_{1\leq i,j <+\infty}E_{X} |h(X_{i},X_{j})|<+\infty.
\end{equation}
(a) For arbitrary $\delta>0$, when the original sample size $n$ is fixed and the bootstrap sample size $m\rightarrow+\infty$, then

\begin{equation}\label{eq 2}
P_{X|w} \Big( \big| U^{*}_{n,m}-\frac{n-1}{n} \ U_n  \big|>\delta   \Big)\rightarrow 0\ a.s.-P_w.
\end{equation}
\begin{equation}\label{eq 2^'}
P_{X|w} \Big( \big| V^{*}_{n,m}- \ V_n  \big|>\delta   \Big)\rightarrow 0\ in\ probability-P_w.
\end{equation}
Consequently, for fixed $n$, as $m\rightarrow+\infty$
\begin{equation}\label{eq 3}
U^{*}_{n,m}\rightarrow \frac{n-1}{n} \ U_n \ in\ probability-P_{X,w}.
\end{equation}
\begin{equation}\label{eq 3^'}
V^{*}_{n,m}\rightarrow V_n \ in\ probability-P_{X,w}.
\end{equation}
(b) For arbitrary $\delta>0$, when $m \rightarrow +\infty$ as $n\rightarrow +\infty$ in such a way that $m\big/ (n \sqrt{2  \log n}) \rightarrow +\infty$, then both (\ref{eq 2}) and (\ref{eq 2^'})  hold true $a.s.-P_w$ and, consequently,
\begin{equation}\label{eq 4}
U^{*}_{n,m}-U_{n}\rightarrow 0\ in\ probability-P_{X,w}
\end{equation}
\begin{equation}\label{eq 4^'}
V^{*}_{n,m}-V_{n}\rightarrow 0\ in\ probability-P_{X,w}.
\end{equation}

\end{thm}
Part (a) of the preceding theorem is closer in spirit to what bootstrap is expected to do and better suites  the inference one would wish to  make as a result of re-sampling repeatedly from a finite original sample.
\par
The next two results are to demonstrate the application of Theorem  \ref{slln for u-stat fix n} for two types of observations. We first consider i.i.d. observations and in Corollary \ref{Corollary u-stat iid} we establish the validity of bootstrap $u$-statistic $U^{*}_{n,m}$. Then, in Corollary \ref{Corollary u-stat beta-mixing}, we establish the validity of the method of $m$-out of-$n$ bootstrap for $U^{*}_{n,m}$
when the observations posses the property of absolute regularity.

\begin{corollary}\label{Corollary u-stat iid}
Let $X_{1}, \ldots,X_n$, which are  the first $n\geq 1$ terms of an infinite  sequence of i.i.d. random variables, and the kernel $h$ be such that $E_{X}|h(X_{1},X_{2})|<+\infty$ and $E_{X}|h(X_{1},X_{1})|<+\infty$. We conclude (a) and (b) as follows.
\\
(a) For arbitrary $\delta>0$, when  $n$ is fixed, as $m\rightarrow +\infty$, (\ref{eq 2}) and (\ref{eq 2^'}) and their respective  consequences (\ref{eq 3}) and (\ref{eq 3^'}) hold true.
\\
(b) For arbitrary $\delta>0$, when $n,m \rightarrow +\infty$ in such a way that $m\big/ (n \sqrt{2  \log n}) \rightarrow +\infty$, then
\begin{equation}\label{eq 5}
P_{X|w} \Big(  \big| U^{*}_{n,m}-E_X h(X_1,X_2)  \big|>\delta   \big)\rightarrow 0 \ a.s.-P_w,
\end{equation}
\begin{equation}\label{eq 6}
P_{X|w}\Big(  \big| V^{*}_{n,m}-E_X h(X_1,X_2)  \big|>\delta   \big)\rightarrow 0 \ a.s.-P_w.
\end{equation}
Consequently,
\begin{equation}\label{eq 5^'}
U^{*}_{n,m}\rightarrow E_X h(X_1,X_2)\ in\ probability-P_{X,w},
\end{equation}

\begin{equation}\label{eq 6^'}
V^{*}_{n,m}\rightarrow E_X h(X_1,X_2)\ in\ probability-P_{X,w}.
\end{equation}

\end{corollary}
\par
It is noteworthy that the diagonal terms $h(X_i,X_i)$, $1\leq i <+\infty $,  have  no influence on  the limit in (\ref{eq 6}) and (\ref{eq 6^'}).
\par
As another example, we  mention  that Theorem \ref{slln for u-stat fix n} can be used to derive the validity of the bootstrap version of the $v$-statistic $V_n$, using the $m$-out of-$n$ scheme. For its use in the next result, we now  consider observations with absolute regularity property. A strong law of large numbers for $u$-statistics based on this type of observations was established by Arcones \cite{Arconese}. In fact part (b) of the following Corollary \ref{Corollary u-stat beta-mixing}  is the bootstrap version of part (i) of Theorem 1 of the latter paper.
\par
To state our next result,  we first   define the mixing coefficient $\beta(k)$ as follows. Consider a strictly stationary sequence of random variables $X_1,X_2,\ldots$,   and  the two $\sigma$-fields $\mathfrak{T}_{1}^{s}:= \sigma(X_1,\ldots,X_s)$ and $\mathfrak{L}_{s}^{\infty}:= \sigma(X_s, X_{s+1},\ldots)$. Let $\{A_i\}_{i=1}^I$ be a partition in $\mathfrak{T}_{1}^{s}$ and $\{B_j\}_{j=1}^J$  be a partition in $\mathfrak{L}_{s+k}^{\infty}$, and define
\begin{equation}\nonumber
\beta(k):= \frac{1}{2}\ \sup\{  \sum_{i=1}^{I} \sum_{j=1}^{J} \big| P_{X}(A_i\cap B_j)-P_X (A_i)P_X(B_j)   \big|; \{A_i\}_{i=1}^I,\ \{B_j\}_{j=1}^J \}.
\end{equation}

\begin{corollary}\label{Corollary u-stat beta-mixing}
Let $X_1,\ldots,X_n$  be the first $n\geq 1$ terms of an infinite sequence of   strictly stationary  random variables.
\\
(a) If the kernel $h$ is  so that $\sup_{1\leq i,j<+\infty} E_{X} |h(X_i,X_j)|<+\infty$, then, for fixed $n$ and arbitrary $\delta>0$, as $m\rightarrow +\infty$, (\ref{eq 2}) and (\ref{eq 2^'}) and their respective  consequences (\ref{eq 3}) and (\ref{eq 3^'}) hold true.
\\
(b) If  for some $s>2$, $\sup_{1\leq i,j<+\infty} E_{X} |h(X_i,X_j)|^{s}<+\infty$, and $\beta(n)\rightarrow 0$, as $n\rightarrow +\infty$, then, for arbitrary $\delta>0 $,  as $n,m\rightarrow +\infty$ in such a way that
$m\big/ (n \sqrt{2  \log n}) \rightarrow +\infty$,
\begin{equation} \label{eq 7}
P_{X|w} \Big( \Big|  U^{*}_{n,m} -\frac{ \sum_{1\leq i \neq j \leq n} E_{X}h(X_i,X_j)}{n(n-1)}    \Big|>\delta    \Big) \rightarrow 0\ a.s.-P_w,
\end{equation}

\begin{equation} \label{eq 8}
P_{X|w} \Big( \Big|  V^{*}_{n,m} -\frac{ \sum_{1\leq i \neq j \leq n} E_{X}h(X_i,X_j)}{n(n-1)}    \Big|>\delta    \Big) \rightarrow 0\ a.s.-P_w.
\end{equation}
Consequently,
\begin{equation}\label{eq 7^'}
U^{*}_{n,m} -\frac{ \sum_{1\leq i \neq j \leq n} E_{X}h(X_i,X_j)}{n(n-1)}\rightarrow 0\ in \ probability-P_{X,w},
\end{equation}

\begin{equation}\label{eq 8^'}
V^{*}_{n,m} -\frac{ \sum_{1\leq i \neq j \leq n} E_{X}h(X_i,X_j)}{n(n-1)}\rightarrow 0\ in \ probability-P_{X,w}.
\end{equation}

\end{corollary}

\section{Asymptotic normality of bootstrap $u$-statistics}\label{section 4}
Conditional (given the observations) central limit theorems for bootstrap  $u$-statistics was  investigated by Arcones and Gin\'{e}  \cite{Arcones Gine}, and  Helmers \cite{Helmers}. However, the closest result in nature to the results in this section is the paper by   Wang and Jing \cite{Wang Jing}, where they provide Edgeworth expansions for the  conditional (given the observations)   distribution  of weighted bootstrap $u$-statistics based on i.i.d. observations,  assuming that the third moment of the kernel $h$ exists and is finite. In contrast, on assuming that the kernel $h$ has a finite second moment, we employ the method of conditioning  on the bootstrap weights, which was first introduced by Cs\"{o}rg\H{o} \emph{et al.}  \cite{Another look}, to  derive a conditional limit theorem for both bootstrap $u$ and $v$-statistics.
\par
Unlike the case of consistency of $u$-statistics, the concept of degeneracy is quite important in studying the asymptotic distribution of  $u$-statistics and their associated $v$-statistics. This is a result of employing the celebrated method of Hoeffding reduction (decomposition) that involves the use of the projections of the underlying $u$-statistic. The number of non-zero projections determine the degree of degeneracy of a $u$-statistic (cf. for example  Serfling \cite{Serfling} or Borovskikh \cite{borovskikh}).   In our case we shall consider the  $u$- and $v$-statistics $U_n$ and $V_n$ defined in (\ref{definition of u-satat}) and (\ref{definition of v-satat}) which are non-degenerate, i.e., for all $i\geq 1$,
$$
\tilde{h}(X_i):=E_{X} \big(   h(X_i,X_j)-E_{X}h(X_i,X_j) \big| X_i  \big)\neq 0\ a.s.-P_X,
 $$
where  $  1\leq j \neq i$.
\par
The concept of degeneracy  is inherited  by the $u$-statistic $U_n$ from its kernel $h$. By this, and  in a similar fashion, the concept of degeneracy can be extended to be used for  \emph{weighted $u$-statistics} as in Nasari \cite{SLLN weighted u-stat}. Degenerate weighted $u$-statistics are to be used in our proofs.
\par
Let $T^{*}_{n,m}$ be the bootstrap version of an original statistic $T_n$, when the $m$-out of-$n$ scheme of bootstrap is used.   The validity in distribution of the $m$-out of-$n$ method of bootstrap should  \emph{ideally} be investigated  by directly showing that, for fixed sample size $n$ and all $t\in \mathds{R}$, as only $m\rightarrow +\infty$, one has
\begin{equation}\label{ideal bootstrap}
 P_{X,w}(T^{*}_{n,m}\leq t) \rightarrow P_{X}(T_n\leq t).
\end{equation}
This  is of the same spirit as part (a) of Theorem \ref{slln for u-stat fix n},   and its Corollaries \ref{Corollary u-stat iid} and \ref{Corollary u-stat beta-mixing} in the previous section. However, unlike the latter results, (\ref{ideal bootstrap}) cannot  be directly investigated as such investigation would require   knowledge of  the sampling distribution of the statistic $T_n$ for fixed $n\geq 1$, which usually   is not the case,  as in practice the sampling distribution is unknown. Therefore,   (\ref{ideal bootstrap}) is   usually  established  by  showing the nearness of $P_{X,w}(T^{*}_{n,m}\leq t)$ to the limiting distribution of the original statistic $T_n$, as $m,n\rightarrow +\infty$. Denoting  the limiting distribution of properly normalized and (usually) centered $T_n $ by $F$, in the  literature,    as  $n,m\rightarrow +\infty$, the stronger conditional (given the sample) version
\begin{equation}\label{ideal bootstrap 2}
P(b_n(T^{*}_{n,m}-T_{n}) \leq t\big| X_1,\ldots,X_n)\rightarrow F(t)\ a.s.-P_X \ or \ in \ probability-P_X
 \end{equation}
is what is shown to  establish  (\ref{ideal bootstrap}),  where $b_n$ is a normalizing sequence. For the Student $t$-statistic, when the limiting distribution $F$ is standard normal, Cs\"{o}rg\H{o} \emph{et al}. \cite{Another look} established the validity of $m$-out of-$n$ scheme of bootstrap via the classical method of (\ref{ideal bootstrap 2}),  as well as using a new method of conditioning on the weights, by showing that, as $n,m\rightarrow +\infty$,
\begin{equation}\nonumber
P_{X|w} (T^{*}_{n,m} \leq t)\rightarrow F(t) \ in\ probability-P_w.
\end{equation}
In this section, for i.i.d. observations, we shall establish the validity, in distribution, of the bootstrap $u$- and $v$-statistics by proving a similar result to the preceding relation, i.e., via conditioning on the bootstrap weights.

\begin{thm}\label{CLT bootstrap u-stat}
Let $X_1,\ldots,X_n$ be the first $n\geq 1$ terms of an infinite  sequence of i.i.d. random variables. Assume that $U_n$ and $V_n$ are non-degenerate.
\\
(a) If $E_{X}h^2(X_1,X_2) <+\infty $, then, for all $t\in \mathds{R}$, as $m,n\rightarrow +\infty$ in such a way that $m=o(n^2)$,
\begin{equation}\label{eq 9}
P_{X|w}\big(  \frac{U^{*}_{n,m}-U_n}{2 \ \hat{\sigma}_n\ \sqrt{\sum_{t=1}^n \big(  \frac{w^{(n)}_t}{m}-\frac{1}{n} \big)^2  }}\leq t  \big) \rightarrow \Phi(t)\ in \ probability-P_w
\end{equation}
(b) If $E_{X}h^2(X_1,X_2),\ E_{X} h^2(X_1,X_1) <+\infty $, then, for all $t\in \mathds{R}$, as $m,n\rightarrow +\infty$ in such a way that $m=o(n^2)$,
\begin{equation}\label{eq 10}
P_{X|w}\big(  \frac{V^{*}_{n,m}-V_n}{2 \ \hat{\sigma}_n\ \sqrt{\sum_{t=1}^n \big(  \frac{w^{(n)}_t}{m}-\frac{1}{n} \big)^2  }}\leq t  \big) \rightarrow \Phi(t)\ in \ probability-P_w,
\end{equation}
where $\Phi$ is the distribution function of  a standard normal random variable,
\begin{equation}\label{sigma hat}
\hat{\sigma}_{n}^2=n (n-1)\sum_{i=1}^{n}(U^{i}_{n-1}-U_{n})^{2},
\end{equation}
and $U^{i}_{n-1}$ is the jackknifed version of $U_{n}$,  based on $X_{1},\ldots,X_{i-1},
X_{i+1},$ $\ldots,X_{n}$, defined  as follows
$$
U^{i}_{n-1}:=\displaystyle{\frac{1}{{n-1 \choose 2}}}
\sum_{\substack{1\leq j_{1}<j_{2}\leq n \\ j_{1}, j_{2}\neq i}} h(X_{j_{1}},X_{j_{2}}).
$$
\end{thm}
\begin{remark}
$\hat{\sigma}_n^{2}$ is the jackknife estimator of  $\sigma^2:=Var(\tilde{h}(X_1))$. It was used by Cs\"{o}rg\H{o} et al.  \cite{Acta mathematica}   for $u$-type processes    when studying the  changepoint  problem  via Studentization.  The estimator $\hat{\sigma}_{n}^2 $
 was generalized by Nasari \cite{Nasari Brazilian} for processes of   $u$-statistics of order greater than or equal to 2. It is noteworthy that $\hat{\sigma}_{n}^2$ remains the right normalizing sequence for the weak convergence (and central limit theorem as a result) of $u$-statistics even when $Eh^{2}(X_1,X_2)=+\infty$ (cf. Theorem 4 and its Corollary 1 of Nasari \cite{Nasari Brazilian}).
\end{remark}
We call attention to  Theorem \ref{CLT bootstrap u-stat} being valid when  bootstrap sub-samples  of size $m=n$, or smaller than $n$,  are  drawn from the original sample.

\begin{remark}
In view of Lemma 1.2 in S. Cs\"{o}rg\H{o} and Rosalsky  \cite{Csorgo Rosalsky}, we note that the respective  unconditional versions of (\ref{eq 9}) and  (\ref{eq 10}) continue to hold true, i.e., having $P_{X,w}$ instead of $P_{X|w}$ in their respective statements with $\Phi$ as their limiting distribution.
\end{remark}

\section{ Remarks on extending the results to $p$-dimensional case} \label{section general case}
  Theorems \ref{SLLN for arrays},  \ref{SLLN for arrays fixed n} and \ref{revers martingale} continue to    hold true for sums of randomly weighted   $p$-dimensional  arrays of randomly weighted  random variables, where $p$ is a positive integer such that $p\geq 2$. Theorems \ref{Marcinkiewicz}, \ref{slln for u-stat fix n} and Corollaries \ref{Corollary u-stat iid} and
  \ref{Corollary u-stat beta-mixing} also hold true for $u-$ and $v-$statistics of order $p\geq 2$.  Mutatis mutandis,  the respective proofs  of these results in the  $p$-dimensional case  remain the same as in the 2-dimensional case of the present paper.  As for Theorem \ref{CLT bootstrap u-stat}, when dealing with $u$-statistics of order $p\geq 2$, the  theorem remains valid assuming the same conditions. The  only change required is to change the constant 2 in the denumerators  of (\ref{eq 9}) and (\ref{eq 10}) to $p$. The proof of this theorem for non-degenerate  $u$, $v$-statistics of order $p$ is in principle   the same (cf. Appendix).

\section{Proofs}\label{section 5}

\subsection{Proof of Theorem \ref{SLLN for arrays}}

The proof of this theorem can be done by modifying the proof of Theorem 1 of \cite{rosalsky} as follows:
\begin{eqnarray*}
&|&\hspace{-0.2 cm} \sum_{1\leq i,j \leq n}  \varepsilon_{ij}^{(n,m)}  X_{ij}^{(n,m)}- \sum_{1\leq i,j \leq n}  c_{ij}^{(n,m)}  X_{ij}^{(n,m)}  |\\
 &\leq&  \sum_{1\leq i,j \leq n} |\varepsilon_{ij}^{(n,m)}- c_{ij}^{(n,m)}   | \ |X_{ij}^{(n,m)}| \\
&\leq& \delta \sum_{1\leq i,j \leq n}  a_{ij}^{(n,m)} |X_{ij}^{(n,m)}|\rightarrow 0\ a.s.-P_{X},
\end{eqnarray*}
as $n\rightarrow \infty$, $\delta\rightarrow 0$. $\square$

\subsection{Proof of Theorem \ref{SLLN for arrays fixed n} }
The proof of this theorem is relatively  simple and we give the details  only for the case when the convergence is a.s.$-P_\varepsilon$.
\par
First consider the case when  $n$ is fixed and $m\rightarrow +\infty$, i.e.,  case (i), and, for $\delta_1,\delta_2>0$, write
\begin{eqnarray*}
&&P_{w} \big\{ \cap_{m\geq 1} \cup_{k\geq m}  P_{X|\varepsilon} \big(   \big|   \sum_{1\leq i,j \leq n}  \varepsilon_{i,j}^{(n,k)}  X_{i,j}^{(n,k)}   \big|>\delta_1      \big)        >\delta_2     \big\}\\
&\leq& P_{w} \big\{ \cap_{m\geq 1} \cup_{k\geq m} \sum_{1\leq i,j \leq n} \big| \varepsilon_{i,j}^{(n,k)} \big|      E_{X} \big|X_{i,j}^{(n,k)}  \big| > \delta_1 \delta_2     \big\}\\
&\leq& P_{w} \big\{ \cap_{m\geq 1} \cup_{k\geq m} \sum_{1\leq i,j \leq n} \big| \varepsilon_{i,j}^{(n,k)} \big|      c_{i,j}^{(n,k)}   > \delta_1 \delta_2     \big\}\\
&=& 0.
\end{eqnarray*}
The latter conclusion is true in view of  assumption (iv), when (i) holds.
\par
On exchanging $m$ by $n$ in the preceding argument, the proof of this theorem when (ii) holds will follow.  Hence the details are omitted.
\par
When (iii) holds,  $m=m(n)\rightarrow+\infty$, as $n$ does so, and with  $\delta_1,\delta_2>0$,  similarly  to the proof for the case (i),  we write
\begin{eqnarray*}
&&P_{w} \big\{ \cap_{n\geq 1} \cup_{k\geq n}  P_{X|\varepsilon} \big(   \big|   \sum_{1\leq i,j \leq k}  \varepsilon_{i,j}^{(k,m(k))}  X_{i,j}^{(k,m(k))}   \big|>\delta_1      \big)        >\delta_2     \big\}\\
&\leq& P_{w} \big\{ \cap_{n\geq 1} \cup_{k\geq n} \sum_{1\leq i,j \leq k} \big| \varepsilon_{i,j}^{(k,m(k))} \big|      E_{X} \big|X_{i,j}^{(k,m(k))}  \big| > \delta_1 \delta_2     \big\}\\
&\leq& P_{w} \big\{ \cap_{n\geq 1} \cup_{k\geq n} \sum_{1\leq i,j \leq k} \big| \varepsilon_{i,j}^{(k,m(k))} \big|      c_{i,j}^{(k,m(k))}   > \delta_1 \delta_2     \big\}\\
&=& 0.
\end{eqnarray*}
Once again, the latter conclusion  is true in view of  assumption (iv), when (iii) holds. This also completes  the proof of Theorem \ref{SLLN for arrays fixed n}. $\square$

\subsection{Proof of Theorem \ref{revers martingale}}
The proof of this theorem  relies on the concept of symmetric (permutable) events, and is similar to the proof of the SLLN for $U$-statistics due Berk \cite{Berk}. To prove this theorem, we first define the $\sigma$-field $\mathfrak{F}_{n}$ as follows:
$$\mathfrak{F}_{n}=\sigma(S_{n},S_{n+1},\ldots).$$
Note that, by definition, $\mathfrak{F}_{n+1}\subset\mathfrak{F}_{n}$, and observe that, since  the $\epsilon_{ij} X_{ij}$'s  are identically distributed, for each $n$ and  each $1\leq  i\neq j\leq n$ we have
\begin{eqnarray*}
E(\epsilon_{12} X_{12}| \mathfrak{F}_{n})&=& E (\epsilon_{ij} X_{ij} | \mathfrak{F}_{n})\ a.s.- P_{X,\epsilon}\\
&=& E(\frac{1}{n(n-1)} \sum_{1\leq i\neq j\leq n} \epsilon_{ij} X_{ij} | \mathfrak{F}_{n} )\ a.s.- P_{X,\epsilon}\\
&=& E(S_{n} | \mathfrak{F}_{n})\ a.s.- P_{X,\epsilon}\\
&=& S_{n} \ a.s.- P_{X,\epsilon}.
\end{eqnarray*}
In view of the latter,  and by Theorem 2.8.6 of \cite{stout}, we now conclude that, as  $n\rightarrow \infty$, we have
$$ S_{n}=E(\epsilon_{12} X_{12}| \mathfrak{F}_{n})\longrightarrow E(\epsilon_{12} X_{12}| \cap_{2\leq n < \infty} \mathfrak{F}_{n})\ a.s.-P_{X,\epsilon}. $$
Now, due to symmetry (permutability) of the tail events (cf.  \cite{stout}),  it follows by the Hewitt-Savage 0-1 law (cf. for example Stout \cite{stout} and Theorem 2.12.4 therein)  that   $E(\epsilon_{12} X_{12}| \cap_{2\leq n < \infty} \mathfrak{F}_{n})$ is a constant a.s.-$P_{X,\epsilon}$. This, in turn  means that, as $n\rightarrow +\infty$,
$$S_{n}\rightarrow E_{X,\epsilon}(\epsilon_{12}X_{12})\ a.s.- P_{X,\epsilon}. $$
The proof of Theorem \ref{revers martingale} is now complete. $\square$

\subsection{Proof of Theorem \ref{Marcinkiewicz}}
To prove this result,  by virtue of (\ref{relation u- and v})  we first observe that
\begin{eqnarray}
m^{-d+2} U^{*}_{n,m}&=& \sum_{1\leq i \neq j \leq n} \frac{w^{(n)}_i w^{(n)}_j}{m^{d-1}(m-1)} h(X_i,X_j)\nonumber\\
&+&\sum_{i=1}^n \frac{w^{(n)^2}_i }{m^{d-1}(m-1)} h(X_i,X_i).\label{eq some number}
\end{eqnarray}
Also, from Theorem 1 of Gin\'{e} and Zinn \cite{Gine Zinn}, we know that, as $n\rightarrow +\infty$,  $n^{-d}\sum_{1\leq i \neq j \leq n} h(X_i,X_j)\rightarrow 0$. Also, observe that, for $\delta>0$,
\begin{eqnarray*}
P_{X}\Big( \frac{\big| \sum_{i=1}^n h(X_i,X_i)  \big|}{n^{d} \log^{2} n} > \delta \Big)
&\leq& \delta^{-2/d} \frac{n}{n^2 \log^2 n} E_X \Big| h(X_1,X_1) \Big|^{2/d} \\
&&=  \delta^{-2/d} \frac{E_X \Big| h(X_1,X_1) \Big|^{2/d}}{n \log^2 n}.
\end{eqnarray*}
Now the proof of Theorem \ref{Marcinkiewicz} results from an application of  Theorem \ref{SLLN for arrays} to both terms on the R.H.S. of the equality (\ref{eq some number}), on taking\\
 $c^{(n,m)}_{ij}=E_{w}\Big(  \frac{w^{(n)}_i w^{(n)}_j}{m^{d-1}(m-1)}\Big)=$ $\left\{
  \begin{array}{ll}
    n^{-2} m^{-d+2}, & \hbox{$i\neq j$;} \\
    \frac{1-1/n}{n m^{d-2}(m-1)}+\frac{1}{n^{2} m^{d-3}(m-1)}, & \hbox{$i=j$}
  \end{array}
\right.
$
\\ and \\
$a^{(n,m)}_{ij}=$ $\left\{
                       \begin{array}{ll}
                         n^{-d}, & \hbox{$i\neq j$;} \\
                         n^{-d} \log^{-d} n, & \hbox{$i=j$}.
                       \end{array}
                     \right.$
\\
\par
Part (a) follows if,  as $n,m \rightarrow +\infty$ such that $m=O(n^{d/(d-2) }\log^{\frac{d+2}{d-2}} n)$, for $\delta>0$, we show that
\begin{equation}\label{eq some number 1}
\sum_{n\geq 1}n(n-1) P_w \Big( \Big| \frac{w^{(n)}_1 w^{(n)}_2}{m^{d-1}(m-1)}-\frac{1}{n^{2}m^{d-2}} \Big|>\delta n^{-d}  \Big)<+\infty.
\end{equation}
In view  of (\ref{eq some number}),  part (b) will follow if, in addition to (\ref{eq some number 1}),  we prove that
\begin{equation}\label{eq some number 2}
\sum_{n\geq 1}n P_w\Big( \Big|   \frac{(w^{(n)}_1)^2 }{m^{d-1}(m-1)}  -\big( \frac{1-1/n}{n m^{d-2}(m-1)}+\frac{1}{n^{2} m^{d-3}(m-1)} \big) \Big|>\frac{\delta}{n^{d}\log^d n} \Big) < +\infty.
\end{equation}
To establish (\ref{eq some number 1}), we first note that without, loss of generality, we may  assume that, for each $1\leq i\neq j \leq n$ and $n\geq 1$,   $w^{(n)}_i w^{(n)}_j \substack{d\\=}$\ Binomial$\big(m(m-1),1/n^2\big)$. This, in turn, enables us to use Bernstein's inequality to estimate
\begin{eqnarray*}
&&n(n-1) P_w \Big( \Big| \frac{w^{(n)}_1 w^{(n)}_2}{m^{d-1}(m-1)}-\frac{1}{n^{2}m^{d-2}} \Big|>\delta n^{-d}  \Big)\\
&\leq&n(n-1) \exp\{ - m(m-1) \frac{\delta^2 m^{2d-4}/ n^{2d}}{2(\frac{1}{n^2}+ \delta m^{d-2}/n^{d})}   \} \\
& \sim& n^2 \exp\{ - \ \frac{m^{d}}{n^d} \ \frac{\delta^2}{2 (n^{d-2}/m^{d-2} + \delta )}   \}.
\end{eqnarray*}
The last term above is a general term of a convergent series when   $n,m \rightarrow +\infty$ such that $m=O(n^{d/(d-2) }\log^{\frac{d+2}{d-2}} n)$. This completes the proof of (\ref{eq some number 1}) and hence that of part (a).
\par
To prove (\ref{eq some number 2}) we use Markov's inequality to estimate
\begin{eqnarray*}
&&n P_w\Big( \Big|   \frac{(w^{(n)}_1)^2 }{m^{d-1}(m-1)}  -\big( \frac{1-1/n}{n m^{d-2}(m-1)}+\frac{1}{n^{2} m^{d-3}(m-1)} \big) \Big|>\frac{\delta}{n^{d}\log^d n} \Big)\\
& \leq& 2 \delta^{-1} n^{d+1} \log^d n \big( \frac{1-1/n}{n m^{d-2}(m-1)}+\frac{1}{n^{2} m^{d-3}(m-1)} \big)\\
&\sim& \frac{2n^d \log^d n}{m^{d-1}}+\frac{2 n^{d-1} \log^d n}{m^{d-2}}.
\end{eqnarray*}
The last term above is a general term of a convergent series when   $n,m \rightarrow +\infty$ such that $m=O(n^{d/(d-2) }\log^{\frac{d+2}{d-2}} n)$. Consequently, we conclude   (\ref{eq some number 2}) and, in view of having also (\ref{eq some number 1}), the proof   of part (b) is  complete. This also completes  the proof of Theorem  \ref{Marcinkiewicz}. $\square$

\subsection{Proof of Theorem \ref{slln for u-stat fix n}}
On taking $c^{(n,m)}_{ij}=c_{ij}=E_{X}|h(X_i,X_j)|$ and noting that here $\sup_{1\leq i,j < +\infty}c_{ij}$ $<+\infty$, the proof of (\ref{eq 2}) will  follow from  Theorem \ref{SLLN for arrays fixed n}    and Remark \ref{new 1}  if we show that, as $m\rightarrow +\infty$,
\begin{equation}\label{eq proof 1}
\sum_{1\leq i\leq j\leq n} \big| \frac{w^{(n)}_i w^{(n)}_{j}  }{m(m-1)}-\frac{1}{n^2}  \big|=o(1) \ a.s.-P_{w}.
\end{equation}
In view of Remark \ref{relation between bootstrap u,v}, as $m\rightarrow +\infty$, (\ref{eq 2^'}) results from the preceding (\ref{eq proof 1}) and the following (\ref{eq proof 2}):
\begin{equation}\label{eq proof 2}
\sum_{i=1}^{n} \big| \frac{(w_{i}^{(n)})^2}{m(m-1)}-\frac{1}{n^2}  \big|=o_{P_{w}}(1).
\end{equation}
To prove (\ref{eq proof 1}), we use Bernstien's inequality along the lines of the following argument to write
\begin{eqnarray}
&&P_{w}\Big(   \sum_{1\leq i \neq j \leq n}  \big| \frac{w^{(n)}_i w^{(n)}_j}{m(m-1)}- \frac{1}{n^2} \big|>\delta  \Big)\nonumber\\
&&\leq n(n-1) P_{w}\Big( \big|    \frac{w^{(n)}_1 w^{(n)}_2}{m(m-1)}- \frac{1}{n^2} \big|>\frac{\delta}{n(n-1)}  \Big)\nonumber\\
&&\leq n(n-1) \exp\big\{ \frac{- m(m-1)}{n(n-1)}. \frac{\delta^2}{2(\frac{n(n-1)}{n^2}+\delta)}   \big\}.\label{eq proof 3}
\end{eqnarray}
Observe now that for the latter upper bound in  (\ref{eq proof 3}), as $n$ is fixed, we have
$$\sum_{m\geq 1} \exp\big\{ \frac{- m(m-1)}{n(n-1)}. \frac{\delta^2}{2(\frac{n(n-1)}{n^2}+\delta)} \big\}<+\infty.$$
 This completes  the proof of (\ref{eq proof 1}).
\par
To establish (\ref{eq proof 2}), we begin with an application of Chebyshev inequality to write
\begin{eqnarray}
P_{w} \Big( \sum_{i=1}^{n} \big| \frac{(w^{(n)}_i)^2}{m(m-1)}-\frac{1}{n^2}  \big|   >\delta  \Big)
& \leq& \delta^{-2} n^2 E_{w} \big(  \frac{(w^{(n)}_i)^2}{m(m-1)}-\frac{1}{n^2} \big)^{2}\nonumber\\
& \sim& \delta^{-2} n^2 \Big( \frac{1}{n m^3}+\frac{7}{n^2 m^2}+\frac{6}{m n^3} -\frac{2(1-\frac{1}{n})}{n^3 m}  \Big)\nonumber\\
&\leq& \delta^{-2} n^2 \Big( \frac{1}{n m^3}+\frac{7}{n^2 m^2}+\frac{6}{m n^3}  \Big).\label{eq proof 4}
\end{eqnarray}
The preceding relation approaches zero as $m\rightarrow +\infty$. This completes the proof of (\ref{eq proof 2}).  Noting now that, in view of Lemma 1.2 in S. Cs\"{o}rg\H{o} and Rosalsky \cite{Csorgo Rosalsky},   (\ref{eq 3})  and  (\ref{eq 3^'}) follow  from (\ref{eq 2}) and (\ref{eq 2^'}), respectively, as $m\rightarrow +\infty$,  completes the proof  of part (a) of Theorem \ref{slln for u-stat fix n}.
\par
To prove part (b) of this theorem, it suffices to observe that the respective upper bounds in   (\ref{eq proof 3}) and (\ref{eq proof 4})  are both general terms of a convergent series in $n$, when the bootstrap sample size $m$ is so that  $m/(n\sqrt{2 \log n})\rightarrow +\infty$. In other words,
\begin{equation}\nonumber
\sum_{n\geq 1} n(n-1) \exp\big\{ \frac{- m(m-1)}{n(n-1)}. \frac{\delta^2}{2(\frac{n(n-1)}{n^2}+\delta)}   \big\}<+\infty
\end{equation}
and
\begin{equation*}
\sum_{n\geq 1} n^2 \Big( \frac{1}{n m^3}+\frac{7}{n^2 m^2}+\frac{6}{m n^3}  \Big)<+\infty.
\end{equation*}
Once again,  the relations (\ref{eq 3}) and (\ref{eq 3^'}) result from Lemma 1.2 in S. Cs\"{o}rg\H{o} and Rosalsky \cite{Csorgo Rosalsky},  and this completes the proof of part (b) and that of Theorem \ref{slln for u-stat fix n}. $\square$

\subsection{Proof of Corollary \ref{Corollary u-stat iid} }
Part (a) is a trivial consequence to Theorem \ref{slln for u-stat fix n}. Here we only have to argue part(b).
\par
From the strong law of large numbers for $u$-statistics based on i.i.d. observations we know that, as $n\rightarrow+\infty$, $U_n\rightarrow E_{X}h(X_1,X_2)$ a.s.-$P_{X}$ (cf. for example Serfling \cite{Serfling}). This combined with (\ref{eq 2}) imply  (\ref{eq 5}) and its consequence (\ref{eq 5^'}).
\par
By virtue of (\ref{relation u- and v}), the proof of (\ref{eq 6}) and its consequence (\ref{eq 6^'}) follow from the fact that $X_i$s are i.i.d. random variables and $E_{X}|h(X_1,X_1)|<+\infty$ and as a result, as $n\rightarrow +\infty$,
$$
\frac{1}{n^2} \sum_{i=1}^n h(X_i,X_i) \rightarrow 0\ a.s.-P_X.
$$

\subsection{Proof of Corollary \ref{Corollary u-stat beta-mixing} }
Once again Part (a) is a trivial consequence of Theorem \ref{slln for u-stat fix n}.
\par
From Theorem 1 of Arcones \cite{Arconese} we know that, as $n\rightarrow +\infty$,  $U_n -\frac{1}{n(n-1)} \sum_{1\leq i \neq j \leq n} E_{X}h(X_i,X_j) \rightarrow 0$ a.s.-$P_X$. This combined with (\ref{eq 2}) imply (\ref{eq 7}) and its consequence (\ref{eq 7^'}).
\par
In view of (\ref{relation u- and v}), to prove  (\ref{eq 8}) and its consequence (\ref{eq 8^'}), it suffices to show that
$$
\frac{1}{n^2} \sum_{1\leq i \leq n} h(X_i,X_i)\rightarrow 0\ a.s.-P_X.
$$
To establish the latter, for $\delta>0$, we use the Chebyshev inequality to write

\begin{eqnarray*}
P_{w}\Big( \big| \frac{1}{n^2} \sum_{1\leq i \leq n} h(X_i,X_i)  \big|>\delta  \Big)&\leq& \delta^{-2} n^{-4} \sum_{1\leq i,j\leq n} E_{X}\big| h(X_i,X_i) h(X_j,X_j)\big |\\
&\leq&  \delta^{-2} n^{-4} \sum_{1\leq i,j\leq n} E^{1/2}_{X} h^{2}(X_i,X_i) \ E^{1/2}_{X} h^{2}(X_j,X_j)\\
&\leq& \delta^{-2}\ n^{-2}  \big(   \sup_{1\leq i<+\infty}E^{1/2}_{X} h^{2}(X_i,X_i) \big)^2 .
\end{eqnarray*}
Now the proof of Corollary \ref{Corollary u-stat beta-mixing} is complete. $\square$

\subsection{proof of Theorem \ref{CLT bootstrap u-stat}}
By virtue of Remark 3 of Nasari \cite{Nasari Brazilian} we may replace $\hat{\sigma}_{n}^2$ by $\sigma^2=E_{X}(\tilde{h})^2\leq E_X h^{2}(X_1,X_2)<+\infty$. This is so since, as $n\rightarrow +\infty$, $\hat{\sigma}_{n}^2\rightarrow \sigma^2$ in probability-$P_X$. Therefore, to prove part (a) of this theorem,  it suffices  to show that
\begin{equation}\label{eq proof 5}
P_{X|w}\big(  \frac{U^{*}_{n,m}-U_n}{2 \sigma\ \sqrt{\sum_{t=1}^n \big(  \frac{w^{(n)}_t}{m}-\frac{1}{n} \big)^2  }}\leq t  \big) \rightarrow \Phi(t)\ in \ probability-P_w
\end{equation}
  Observe now that $U_{n,m}^{*}-\ U_n$  can be written as a weighted $u$-statistic as follows:
\begin{equation}\label{eq proof 6}
U_{n,m}^{*}-\ U_n= \sum_{1\leq i \neq j \leq n} \big(\frac{w^{(n)}_i w^{(n)}_j}{m(m-1)}-\frac{1}{n(n-1)}   \big) h(X_i,X_j).
\end{equation}
In view of the preceding relation,  the fact that $\sum_{1\leq i \neq j \leq n} \frac{w^{(n)}_i w^{(n)}_j}{m(m-1)}=1$ allows us to assume  without loss of generality  that $E_{X}h(X_1,X_2)=0$. In other words, the difference $U_{n,m}^{*}-U_n$ does not feel the theoretical mean of the original $u$-statistic $U_n$.
\par
We now employ a Hoeffding type reduction for the weighted $u$-statistic $U_{n,m}^{*}-U_n$, as in (\ref{eq proof 6}), and  write
\begin{eqnarray}
U_{n,m}^{*}-U_n&=& \sum_{1\leq i \neq j \leq n} \big(\frac{w^{(n)}_i w^{(n)}_j}{m(m-1)}-\frac{1}{n(n-1)}   \big) \big(h(X_i,X_j)-\tilde{h}(X_i)-\tilde{h}(X_j)\big)\nonumber\\
&+& 2 \sum_{1\leq i \neq j \leq n} \big(\frac{w^{(n)}_i w^{(n)}_j}{m(m-1)}-\frac{1}{n(n-1)}   \big) \tilde{h}(X_i) \nonumber\\
&=:& G_{n,m}+H_{n,m}.
\end{eqnarray}
\par
The proof of part (a) follows if we show that, as $n,m\rightarrow +\infty$ such that $m=o(n^2)$, then, for $\delta_1,\delta_2>0$,
\begin{equation}
P_w\big\{P_{X|w} (\frac{| G_{n,m} |}{\sqrt{\sum_{t= 1}^n (\frac{w_{t}^{(n)}}{m}-\frac{1}{n}\big)^2 }}>\delta_1)>\delta_2 \big\}=o(1),\label{eq proof 7}
\end{equation}
and for all $t\in \mathbb{R}$,
\begin{equation} P_{X|w} \Big(\frac{H_{n,m}}{\sigma \sqrt{ \sum_{t= 1}^n (\frac{w_{t}^{(n)}}{m}-\frac{1}{n}\big)^2 }}\leq t \Big)\rightarrow \Phi(t) \ in \ probability-P_w.\label{eq proof 9}
\end{equation}
\par
We are now to prove (\ref{eq proof 7}). To do so   we first show that as  $n,m\rightarrow +\infty$ such that $m=o(n^2)$, for $\delta>0$
\begin{equation}\label{eq proof 8}
  P_{w}\big(    \big| \frac{m}{(1-\frac{1}{n})} \sum_{t=1}^{n} \big( \frac{w_t^{(n)}}{m} - \frac{1}{n} \big)^2 -1 \big|> \delta  \big)\rightarrow 0.
\end{equation}
First note that for each $i$,  $1\leq i \leq n$,
$$\displaystyle{ E_w (\frac{w_t^{(n)}}{m} - \frac{1}{n})^{2}= E_w (\frac{w_1^{(n)}}{m} - \frac{1}{n})^{2}=\frac{(1-\frac{1}{n})}{n m}  }.$$
We now employ the Chebyshev  inequality  to bound the $P_{w}(.)$ statement  of (\ref{eq proof 8}) as follows.
\begin{eqnarray*}
&&  P_{w}\big(    \big| \frac{m}{(1-\frac{1}{n})} \sum_{t=1}^{n} \big( \frac{w_t^{(n)}}{m} - \frac{1}{n} \big)^2 -1 \big|> \delta  \big)  \\
&&\leq \frac{m^{2}}{\delta^2  (1-\frac{1}{n})^2 } E_w\Big(\sum_{t=1}^{n} \big( \frac{w_t^{(n)}}{m} - \frac{1}{n})^{2}-\frac{(1-\frac{1}{n})}{m}    \Big)^{2}\\
&=& \frac{m^{2}}{\delta^{2}  (1-\frac{1}{n})^2}  \Big\{ E_w\Big(\sum_{t=1}^{n} \big( \frac{w_{t}^{(n)}}{m} - \frac{1}{n})^{2} \Big)^{2}- \frac{(1-\frac{1}{n})^2}{m^{2}_{n}}   \Big\}\\
&=& \frac{m^{2}_{n}}{\delta^2  (1-\frac{1}{n})^2 } \Big\{ n E_w \big( \frac{w^{(n)}_{1}}{m}-\frac{1}{n} \big)^{4}\\
&& + n(n-1) E_w \Big( \big(\frac{w^{(n)}_{1}}{m}-\frac{1}{n}\big)^2 \big(\frac{w^{(n)}_{2}}{m}-\frac{1}{n}\big)^2    \Big) - \frac{(1-\frac{1}{n})^2}{m^{2}_{n}}   \Big\}.
\end{eqnarray*}
 In view of the fact that  $w^{(n)}_{t}$, $1\leq t \leq n$,  have multinomial distribution,          after computing  $E_{w} \big[(w^{(n)}_1)^{a} (w^{(n)}_2)^{b}\big]$, where $a,b$ are  two integers such that  $0\leq a,b\leq 2$, followed by some algebra,  we can bound the preceding relation above by:
 \begin{eqnarray*}
&&\frac{m^{2}}{\delta^2  (1-\frac{1}{n})^2 } \Big\{  \frac{(1-\frac{1}{n})}{n^{3}m^{3}_{n}}+ \frac{(1-\frac{1}{n})^{4}}{m^{3}_{n}}+\frac{(m -1)(1-\frac{1}{n})^{2}}{n m^{3}}+ \frac{4 (n-1)}{n^{3} m}\\
  && +  \frac{1}{m^{2}}- \frac{1}{n m^{2}}+\frac{n-1}{n^{3}m^{2}} + \frac{4(n-1)}{n^{2} m^{3}}- \frac{(1-\frac{1}{n})^{2}}{m^{2}_n} \Big\}\\
 &&\sim \frac{1}{\delta^2}  \Big\{   \frac{4 m}{n^{2}}+  \frac{1}{n^3 m}+ \frac{1}{m}+  \frac{1}{n^{2}}+ \frac{4}{n m}      \Big\},
\end{eqnarray*}
where $a_n\sim b_n$ stands for the asymptotic equivalence  of the numerical sequences $a_n$ and $b_n$.
\par
Clearly, as $n,m\rightarrow \infty$,  the preceding relation approaches zero when $m=o(n^2)$. Now the  proof of (\ref{eq proof 8}) is complete.
\par
In view of (\ref{eq proof 8}), equivalently to  (\ref{eq proof 7}),  we show that as $n,m \rightarrow +\infty$ so that $m=o(n^2)$, then
\begin{equation}\label{eq proof 7^'}
P_{w} \big\{  P_{X|w}(m^{1/2}   | G_{n,m}| >\delta_1)>\delta_2 \big\}
 \rightarrow 0.
\end{equation}
In order to establish the preceding conclusion,  we first note that $G_{n,m}$ is a weighted $u$-statistic and its kernel $g(X_i,X_j):= h(X_i,X_j)-\tilde{h}(X_i)-\tilde{h}(X_j)$ possesses  the property of (complete) degeneracy (cf. Nasari \cite{SLLN weighted u-stat} and Remark 1 therein). On noting now that part (b) of Lemma 2 of Nasari \cite{SLLN weighted u-stat}  remains true even  when the weights are not necessarily  non-negative, we use it   in the second line of the  following argument to conclude
\begin{eqnarray}
&&P_{w} \big\{  P_{X|w}(m^{1/2}   | G_{n,m}| >\delta_1)>\delta_2 \big\}\nonumber \\
&&\leq P_{w} \big\{  m \sum_{1\leq i \neq j \leq n} \big(\frac{w_{i}^{(n)} w_{j}^{(n)}  }{m(m-1)}-\frac{1}{n(n-1)}   \big)^2> \frac{\delta_1^{2} \delta_2}{A E_{X}h^{2}(X_1,X_2)}   \big\},\label{give it number}
\end{eqnarray}
where $A$ is a positive constant that does not depend on $n$ or $m$. Letting now $\delta^{-1}_3= \frac{A E_{X}h^{2}(X_1,X_2)  }{\delta_1^{2} \delta_2}$, we can bound the R.H.S. of the preceding inequality by
\begin{eqnarray*}
&&\delta^{-1}_3 \ m n(n-1)  \big\{ E_{w} \big( \frac{w_{1}^{(n)} w_{2}^{(n)}}{m(m-1)} \big)^2- \frac{2E_{w}(w_{1}^{(n)} w_{2}^{(n)}) }{nm(n-1)(m-1)}+\frac{1}{n^2(n-1)^2}    \big\}\\
&=&  \delta^{-1}_3 \ m n(n-1)  \big\{  \frac{m(m-1)(m-2)(m-3)}{m^2(m-1)^2 n^{3}(n-1)}+\frac{m(m-1)(m-2)}{m^2(m-1)^2 n^{2}(n-1)}\\
&& +\frac{m(m-1)}{m^2(m-1)^2 n(n-1)}-\frac{2m(m-1)}{m(m-1)n^3 (n-1)}+\frac{1}{n^4} \big\}\\
&\sim& \delta^{-1}_3 (\frac{2}{n}+\frac{1}{m})\rightarrow 0,\ as \ n,m\rightarrow+\infty.
\end{eqnarray*}
Now the proof of (\ref{eq proof 7}) is complete.
\par
To derive (\ref{eq proof 9}), we first write   $H_{n,m}$ as follows.
\begin{eqnarray*}
H_{n,m}&=& \sum_{i=1}^n \tilde{h}(X_i) \big\{ \sum_{j=1}^n \big(\frac{w^{(n)}_i w^{(n)}_j}{m(m-1)}-\frac{1}{n(n-1)}\big)- \big(  \frac{(w_{i}^{(n)})^2}{m(m-1)}-\frac{1}{n(n-1)}\big)  \big\}\\
&=& \sum_{i=1}^n \big(\frac{w_{i}^{(n)}}{m-1}-\frac{1}{n-1}     \big)  \tilde{h}(X_i)+ \sum_{i=1}^n \big(\frac{(w_{i}^{(n)})^2}{m(m-1)}-\frac{1}{n(n-1)}     \big)  \tilde{h}(X_i)\\
&=:& H_{n,m}(1)+H_{n,m}(2).
\end{eqnarray*}
Note that for $H_{n,m}(1)$  we have
\begin{eqnarray*}
H_{n,m}(1)&=& \sum_{i=1}^n \big(\frac{w_{i}^{(n)}}{m}-\frac{1}{n}     \big)  \tilde{h}(X_i)\\
&+& \sum_{i=1}^n \big(\frac{w^{(n)}_i}{m(m-1)}-\frac{1}{n(n-1)} \big) \tilde{h}(X_i).
\end{eqnarray*}
Since, $\tilde{h}(X_1),\tilde{h}(X_2),\ldots$, is a sequence of centered i.i.d. random variables,   from conclusion (2.13)  of Corollary 2.1 of Cs\"{o}rg\H{o} \emph{et al.} \cite{Another look}, on replacing the sample variance $S_n$ by the population variance $\sigma $ therein,  we know that,  as $m,n\rightarrow +\infty$ such that $m=o(n^2)$, then
\begin{equation*}
P_{X|w}\Big(\frac{\sum_{i=1}^n \big(\frac{w_{i}^{(n)}}{m}-\frac{1}{n}     \big)  \tilde{h}(X_i)}{ \sigma \sqrt{\sum_{t=1}^{n} \big( \frac{w_t^{(n)}}{m} - \frac{1}{n} \big)^2 } }\leq t\Big)\rightarrow \Phi(t)\ in \ probability-P_{w}
\end{equation*}
for all $t\in \mathbb{R}$.  By this, in view of (\ref{eq proof 8}), to complete the proof of part (a) it only remains to show that
\begin{equation}
 P_{X|w} ( m^{1/2}  \big|  \sum_{i=1}^n \big(\frac{w_{i}^{(n)}}{m(m-1)}-\frac{1}{n(n-1)} \big) \tilde{h}(X_i) \big|>\delta )=o_{P_w}(1),\label{eq proof 10}
\end{equation}
and
\begin{equation}
 P_{X|w}(m^{1/2} | H_{n,m}(2) |>\delta  )=o_{P_w}(1),\label{eq proof 11}
\end{equation}
where $\delta>0$ is an arbitrary positive constant.

\par
To  establish  relation (\ref{eq proof 10}), for $\delta_,\delta_2>0$,  we write
\begin{eqnarray*}
&&P_{w} \big\{ P_{X|w} \Big( m^{1/2}  \big|  \sum_{i=1}^n \big(\frac{w_{i}^{(n)}}{m(m-1)}-\frac{1}{n(n-1)} \big) \tilde{h}(X_i) \big|>\delta_1 \Big)> \delta_2 \big\}  \\
&\leq & P_{w} \big\{  m^{1/2}\sum_{i=1}^n \big| \frac{w_{i}^{(n)}}{m(m-1)}-\frac{1}{n(n-1)} \big| >   \frac{\delta_1 \delta_2}{2 E_{X}|h(X_1, X_2)|}  \big\}\\
&\leq& \frac{2 E_{X}|h(X_1, X_2)| m^{1/2} n  }{\delta_1 \delta_2}  \big\{   \frac{m}{m(m-1)n}+\frac{1}{n(n-1)} \big\}\\
&\sim& \frac{2 E_{X}|h(X_1, X_2)|}{\delta_1 \delta_2} \big\{ \frac{1}{\sqrt{m}} + \frac{\sqrt{m}}{n}   \big\}
\rightarrow 0,
\end{eqnarray*}
as $m,n \rightarrow +\infty$ such that $m=o(n^2)$.
\par
To prove (\ref{eq proof 11}), with a similar argument as the one used to prove  (\ref{eq proof 10}), we bound its $P_{X|w}(.)$ statement by
\begin{eqnarray*}
&&P_{w} \big\{  m^{1/2} \sum_{i=1}^n \big|  \frac{(w_{i}^{(n)})^2}{m(m-1)}-\frac{1}{n(n-1)}  \big|    >  \frac{\delta_1 \delta_2}{E_{X}|h(X_1,X_2)|} \big\}\\
&\leq& \frac{  {2 E_{X}|h(X_1,X_2)|}\ m^{1/2}\ n }{ {\delta_1 \delta_2} }  \big(  E_{w} ( \frac{(w_{i}^{(n)})^2}{m(m-1)})+\frac{1}{n(n-1)}  \big)\\
&\sim& \frac{  {2 E_{X}|h(X_1,X_2)|} }{ {\delta_1 \delta_2} } \big( \frac{1}{\sqrt{m}}+\frac{2\sqrt{m}}{n} \big) \rightarrow 0.
\end{eqnarray*}
This completes the proof of (\ref{eq proof 11}) and that of part (a) of this theorem.
\par
The proof of part (b)  to a large extent relies on part (a) in view of the following relation.
\begin{eqnarray}
\frac{V^{*}_{n,m}-V_{n}}{2 \hat{\sigma}_n \sqrt{\sum_{t=1}^n ( \frac{w_{t}^{(n)}}{m}-\frac{1}{n}   )^2 } }&=& \frac{U_{n,m}^{*}-U_{n}}{2 \hat{\sigma}_n \sqrt{\sum_{t=1}^n ( \frac{w_{t}^{(n)}}{m}-\frac{1}{n}   )^2 } }\label{eq proof 12} \\
&+& \frac{\sum_{1\leq i\neq j \leq } h(X_i,X_j) }{ 2 \hat{\sigma}_{n}  n^2 (n-1)    \sqrt{\sum_{t=1}^n ( \frac{w_{t}^{(n)}}{m}-\frac{1}{n}   )^2 }    }\label{eq proof 13}\\
&+& \frac{\sum_{i=1}^n \big( \frac{(w_{i}^{(n)})^2}{m(m-1)}-\frac{1}{n^2}\big) h(X_i,X_i)   }{2 \hat{\sigma}_n \sqrt{\sum_{t=1}^n ( \frac{w_{t}^{(n)}}{m}-\frac{1}{n}   )^2 } }.\label{eq proof 14}
\end{eqnarray}
In part (a) it was shown that the  conditional distribution of the (\ref{eq proof 12}) part of the above equality  converges to normal in probability-$P_w$. We now show the  hierarchical asymptotic negligibility of parts  (\ref{eq proof 13}) and (\ref{eq proof 14}). To do so, we first replace the estimator $\hat{\sigma}^{2}_n$ by $\sigma^2=E_{X}\tilde{h}^2(X_1)$ and  $\sum_{t=1}^n ( \frac{w_{t}^{(n)}}{m}-\frac{1}{n}   )^2$ by
$1/m$.
\par
To deal with (\ref{eq proof 13}),  it suffices to observe the asymptotic a.s.-$P_X$ negligibility of
\begin{equation*}
\frac{m^{1/2}  }{n} . \frac{\sum_{1\leq i\neq j \leq n}  h(X_i,X_j)}{n(n-1)}
\end{equation*}
that results from the strong law of large numbers for $u$-statistics (cf. for example Serfling (1980)), provided that $m=o(n^2)$.
\par
We now use a similar  argument to the one used to prove (\ref{eq proof 11}) to write
\begin{eqnarray*}
&& P_w\big\{ P_{X|w}(m^{1/2} | \sum_{i=1}^n \big( \frac{(w_{i}^{(n)})^2}{m(m-1)}-\frac{1}{n^2}\big) h(X_i,X_i) |>\delta_1  )>\delta_2 \}\\
&\leq & P_{w} \big\{  m^{1/2} \sum_{i=1}^n \big|  \frac{(w_{i}^{(n)})^2}{m(m-1)}-\frac{1}{n(n-1)}  \big|    >  \frac{\delta_1 \delta_2}{E_{X}|h(X_1,X_1)|} \big\}\\
&\rightarrow& 0,\ as\ n,m\rightarrow+\infty.
\end{eqnarray*}
The preceding  yields the asymptotic hierarchical negligibility of (\ref{eq proof 14}). By this the proof of part (b) and that of Theorem \ref{CLT bootstrap u-stat} are complete. $\square$

\section*{Appendix}
The use of the Hoeffding decomposition to reduce the underlying weighted $u$-statistic of order $p\geq 3$ to weighted partial sums yields   $p-1$   \emph{completely degenerate}  weighted $u$-statistics of order $p$ to order 2. In view of Part (b) of Lemma 2 of Nasari \cite{SLLN weighted u-stat},  the  hierarchical  asymptotic  negligibility of each one of these completely degenerate weighted $u$-statistics can be done by establishing an  approximation similar to (\ref{give it number}). This should be followed by an application of Chebyshev's inequality, similarly  to the argument that follows (\ref{give it number}). Similarly to the case $p=2$, the conditional asymptotic  normality for  the weighted partial sums of i.i.d. random variables, resulting from the Hoeffding reduction, is concluded  from Corollary 2.1 of Cs\"{o}rg\H{o} \emph{et al}. \cite{Another look} , via  steps similar  to those used when dealing with $H_{n,m}$ in Theorem \ref{CLT bootstrap u-stat}.  $\square$


\bibliographystyle{model1a-num-names}







\end{document}